\documentclass[10pt,a4paper, twoside]{article}
\date{\today}

\usepackage{latexsym,amsmath,amsfonts,amscd, amssymb, mathrsfs, slashed, 
amsthm, epsfig, psfrag}
\usepackage[english]{babel}

\theoremstyle{plain}  % default
\newtheorem{theorem}{Theorem}[section]

\newtheorem*{theorem*}{Theorem}

\newtheorem{corollary}[theorem]{Corollary}
\newtheorem{lemma}[theorem]{Lemma}
\newtheorem{proposition}[theorem]{Proposition}

\newtheorem{definition}[theorem]{Definition}
\theoremstyle{remark}
\newtheorem{example}[theorem]{Example}
\newtheorem{remark}[theorem]{Remark}

\newtheorem*{claim*}{Claim}

\numberwithin{equation}{section}

\renewcommand{\leq}{\leqslant}

\renewcommand{\geq}{\geqslant}

\newcommand{\lra}{\longrightarrow}
\newcommand{\R}{\mathbb{R}}

\newcommand{\C}{\mathbb{C}}

\newcommand{\del}{\partial}
\newcommand{\ds}{\displaystyle}

\def\lms{\longmapsto}
\def\proof{\noindent\textit{Proof ---}\hspace*{0.2cm}}
\def\qed{\vspace*{-0.1cm} \hfill{$\square$}}

\newcommand{\g}{\mathfrak{g}}

\def\Ric{\mathrm{Ric}}

\def\im{\mathrm{im}}
\def\g{\mathfrak{g}}

\def\m{\mathfrak{m}}
\newcommand{\SU}{\mathrm{SU}}
\newcommand{\x}{\times}

\newcommand{\Ad}{\ensuremath{\mathrm{Ad}\,}}
\newcommand{\Spin}{\ensuremath{\mathrm{Spin}}}

\def\haken{\mathbin{\hbox to 6pt{%
                 \vrule height0.4pt width5pt depth0pt
                 \kern-.4pt
                 \vrule height6pt width0.4pt depth0pt\hss}}}

\begin{document}
%%%%%%%%%%%%%%%%%%%%%%%%%%%%%%%%%%%%%%%%%%%%%%%%%%%%%%%%%%%
\title{Einstein manifolds with skew torsion}
%%%%%%%%%%%%%%%%%%%%%%%%%%%%%%%%%%%%%%%%%%%%%%%%%%%%%%%%%%%

\author{Ilka Agricola\footnote{Fachbereich Mathematik und Informatik,
Philipps-Universit\"at Marburg, Campus Hans-Meerwein-Stra\ss e, 35032 Marburg,
Germany, \texttt{agricola@mathematik.uni-marburg.de}} \and  
Ana Cristina Ferreira\footnote{Centro de Matem\'atica, Universidade de
 Minho, Campus de Gualtar, 4710-057 Braga, Portugal, 
\texttt{anaferreira@math.uminho.pt}}}

%\subjclass[2000]{Primary 53C21; Secondary 58A14}

\date{October 1, 2013.}

\maketitle

\begin{abstract}
Abstract. This paper is devoted to the first systematic investigation
of manifolds that are Einstein for a connection $\nabla$ with skew 
symmetric torsion. We derive the 
Einstein equation from a variational principle and
prove that, for parallel torsion, any Einstein manifold with skew torsion
has constant scalar curvature; and if it is
complete of positive scalar $\nabla$-curvature, it is necessarily compact and
it has finite first fundamental group $\pi_1$. 
The longest part of the paper is  devoted to the systematic construction
of large families of examples. We discuss when a Riemannian 
Einstein manifold can be Einstein with skew torsion. We give examples of
almost Hermitian, almost metric contact, and $G_2$ manifolds that are
Einstein with skew torsion. For example, we prove that any Einstein-Sasaki manifold and
any $7$-dimensional 3-Sasakian manifolds admit deformations into
an Einstein metric with parallel skew torsion.

\end{abstract}

%\tableofcontents

%%%%%%%%%%%%%%%%%%%%%%%%%%%%%%%%%%%%%%%%%%%%%%%%%%%%%%%%%%%%%%%%%%%%%%%%%%%
\section{Preliminaries}
%%%%%%%%%%%%%%%%%%%%%%%%%%%%%%%%%%%%%%%%%%%%%%%%%%%%%%%%%%%%%%%%%%%%%%%%%%%

\subsection{Introduction}
%%%%%%%%%%%%%%%%%%%%%%%%%%%%%%%%%%%%%%%%%%%%%%%%%%%%%%%%%%%%%%%%%%%%%%%%%%%%
Torsion, and in particular skew torsion, has been a topic of
interest to both mathematicians and physicists in recent decades.
The first attempts to modify general relativity by introducing torsion
go back to the 1920's with the work of \'{E}. Cartan \cite{Cart23},
and were deepened--in modified form--from the 1970's in Einstein-Cartan theory.
More recently, the torsion of a connection makes its appearance in 
superstring compactifications, 
where the basic model of type II string theory consists of a Riemannian
manifold, a connection with skew torsion, a spinorial field and a
dilaton function.

>From the mathematical point of view, skew torsion has played a
significant role in the proof of the local index theorem
for Hermitian non-K\"{a}hler manifolds \cite{Bism89} and it is
a standard tool for the investigation of non-symmetric homogeneous 
spaces, since the canonical connection of such a space does not
coincide with the Levi-Civita connection anymore \cite{TricVane83}.
In generalized
geometry \cite{Hitc10, Gual03}, there are natural connections
with skew torsion, the exterior derivative of the $B$-field.

Torsion is also ubiquitous in the theory of non-integrable
geometries. This field has been revived in recent years through the
development of superstring theory. Firstly, integrable geometries
(like Calabi-Yau manifolds, Joyce manifolds, etc.) are exact
solutions of the Strominger model with vanishing $B$-field. By
deforming these vacuum equations and looking for models with
non-trivial $B$-field, a new mathematical approach implies that
solutions can be constructed geometrically from non-integrable
geometries with torsion (for example, almost Hermitian, almost
metric contact or weak $G_2$ structures). If $(M,g)$ denotes a 
Riemannian manifold, we will write any metric connection on $M$
as ($\nabla^g$ denotes the Levi-Civita connection)
$$
\nabla_X Y = \nabla_X^g Y + A(X,Y).
$$ 
We say that $\nabla$ has skew torsion if
the  contraction of its torsion $H(X,Y)$  with $g$,  $H(X,Y,Z) := g(H(X,Y),Z)$, 
is totally antisymmetric. In this case, $A(X,Y)= 1/2 H(X,Y,-)$.
>From the three Cartan classes of torsion,
this is by far the richest and the best understood: 
such connections are always complete, they are the only ones with
non-trivial coupling  to the Dirac operator such that the resulting
Dirac operator is still formally self-adjoint \cite{Friedrich79},
and  many non-integrable geometric structures admit a \emph{unique} invariant
connection with totally antisymmetric torsion, thus it is a natural 
replacement for the Levi-Civita connection \cite{FriedIvan02}. 
Also, many new results on the holonomy
properties of connections with skew torsion are now available 
\cite{Olmos&R12}, and a lot of effort has been devoted by many researchers
to the construction of geometrically interesting examples.

\medskip\noindent {\bf Outline.} 
In this article, we propose a notion of `Einstein manifold with skew
torsion' for an $n$-dimensional Riemannian manifold
equipped with a metric connection with skew torsion\footnote{
To prevent any confusion: This is to be understood as a generalization of the
\emph{mathematical} Einstein equation, and not as an alternative
field equation for the gravitational field.}. 
Our approach will be mostly dimension independent and deal, where possible, 
with   general issues; the comparison with the
results obtained previously in General Relativity (Section \ref{Sec:GR})
will illustrate how our approach differs from the previous work
in the area. We start by deducing the Einstein equation with skew torsion
from a variational principle. In order to investigate
the curvature properties of $\nabla$-Einstein manifolds, more assumptions
are needed; for example, easy examples illustrate that, in general, the
scalar curvature will not be constant. We show that a very suitable
restriction is to impose that the torsion be parallel, $\nabla H=0$.
There are several families of manifolds that are  classically known to 
admit parallel characteristic torsion, namely 
nearly K\"ahler manifolds, Sasakian manifolds, nearly parallel $G_2$-manifolds,
and naturally reductive spaces; these classes have been considerably enlarged in
more recent  work (see \cite{Vais79}, \cite{GaudOrne98},
\cite{AlexFrieScho05}, \cite{Alexandrov03}, 
\cite{Frie07}, \cite{Scho07}), leading to a host of instances 
to which our theorems can be applied. The key result  illustrating
that this is the `right' condition is the following: If $\nabla H=0$,
any $\nabla$-Einstein manifold has constant scalar curvature
(both Riemannian and of the connection with torsion); and if it is
complete of positive scalar $\nabla$-curvature, it is necessarily compact and
it has finite first fundamental group $\pi_1$. Thus, we obtain
the best possible analogy to the Riemannian case.
We then discuss an easy, but powerful criterion when a
Riemannian Einstein space will be $\nabla$-Einstein for a given torsion 
$3$-form.

The longest part of the paper is devoted to the systematic construction
of examples in different situations. We first treat the case $n=4$,
where the second author had proposed an alternative definition of 
`Einstein with torsion' based on the phenomenon of self-duality
\cite{Ferr10, Ferr11}. In general, this is a different concept, but
we will show that they coincide if one assumes
parallel torsion. Under this condition, we observe that a $4$-dimensional
Hermitian Einstein manifold is locally isometric to $\R\times S^3$. 
After a quick discussion of the Lie group case, we treat almost
Hermitian manifolds in dimension $6$, where we identify a class of 
homogeneous manifolds of type $W_1\oplus W_3$ that is always
Einstein with parallel skew torsion; this includes, in particular,
all nearly K\"ahler manifolds. For almost contact manifolds,
$\nabla$-Einstein implies $\nabla$-Ricci-flatness, since the contact
distribution is a $\nabla$-parallel vector field. We prove that
\emph{every} Einstein-Sasaki manifold with its characteristic torsion
admits a deformation into a $\nabla$-Einstein-Sasaki manifold. 
Thus, there is a multitude of $\nabla$-Ricci-flat Einstein spaces
in all odd dimensions. The $7$-dimensional case is treated separately
because of its relevance for $G_2$ geometry. Again, all nearly parallel 
$G_2$ manifolds are Einstein with parallel skew torsion; moreover,
we show that any $7$-dimensional $3$-Saskian manifold carries three different
connections that turns it into an Einstein manifold with parallel skew
torsion, and that it admits a deformation of the metric that
carries again an Einstein structure with parallel skew
torsion. Finally, we present several 
$\nabla$-Einstein structures on Aloff-Wallach manifolds $\SU(3)/S^1$;
several of them are new, i.\,e.~not among those that were predicted 
theoretically in the previous sections. 

We end this outline with some conjectural remarks. 
In the past years, there has been a
revived interest in higher dimensional black holes, i.\,e.~Ricci flat
manifolds with Lorentzian signature, because of the exciting discovery
of new horizon topologies (`black rings') and the option to use these
as more sophisticated backgrounds for superstring theories. On the
other side, the use of skew torsion is by now a well-established tool
in superstring compactifications. Thus, we believe that Einstein spaces
with torsion will be of interest for future developments in this area as well,
although the present paper will not deal with these issues.

\medskip\noindent {\bf Acknowledgements.} 
Both authors thank Thomas Friedrich (Berlin) for his steady
mathematical interaction.
Ilka Agricola warmly thanks Friedrich Hehl (K\"oln) for intensive
e-mail discussions and valuable references to the work
done in general relativity on Einstein-Cartan theory.
She acknowledges financial support by the
DFG within the priority programme 1388 "Representation theory".
Ana Ferreira thanks Philipps-Universit\"at Marburg for its
hospitality during a research stay in April-July 2012,
and acknowledges financial support by the DAAD for this stay.
The second author was partially supported by FEDER Funds through Programa 
Operacional Factores de Competitividade --- COMPETE and by Portuguese 
Funds through FCT 
within the Projects PEst-C/MAT/UI0013/2011 and PTDC/MAT/118682/2010.
Finally, we thank Anna Gori (Milano) for pointing out a mistake in a
preliminary version of this article.

%%%%%%%%%%%%%%%%%%%%%%%%%%%%%%%%%%%%%%%%%%%%%%%%%%%%%%%%%%%%%%%%%%%%%%%%%%%%
\subsection{Einstein spaces with torsion in General Relativity}\label{Sec:GR}
%%%%%%%%%%%%%%%%%%%%%%%%%%%%%%%%%%%%%%%%%%%%%%%%%%%%%%%%%%%%%%%%%%%%%%%%%%%%

The first attempts to introduce torsion as an additional 'data' for 
describing physics in general relativity go back to Cartan himself 
\cite{Cartan24a}. Viewing torsion as some intrinsic angular momentum,
 he derived a set of gravitational field equations from a variational
principle, but postulated that the energy-momentum tensor should still
be divergence-free, a condition too  restrictive
for making this approach useful. The idea was taken up again in 
broader context in the late fifties. The variation of the scalar curvature
and of an additional Lagrangian generating the energy-momentum  and
the spin tensors on a 
space-time endowed with a metric connection with torsion
yielded the two fundamental equations of \emph{Einstein-Cartan 
theory},
first formulated by Kibble \cite{Kibble61} and Sciama (see his article 
in \cite{Infeld62}). The first equation can, by some elimination process,
be reduced to an equation which is similar to  Einstein's classical 
field equation of general relativity with an effective energy 
momentum tensor $T_{\mathrm{eff}}$ depending on torsion, the second one  
relates the torsion to the spin density (in the absence of spin, the 
torsion vanishes and the whole theory reduces to Einstein's original 
formulation of general relativity).
A.~Trautman provided an 
elegant formulation of  Einstein-Cartan theory in the
language of principal fibre bundles \cite{Trautman73a}. 
For a general review of gravity with spin and torsion including extensive 
references, we  recommend the article \cite{Hehl&H&K&N76} or the new 
`source' book \cite{BlaHehl12}, which contains most of the articles
cited in this section with extensive commentaries.  

For this article, our main interest will be in general results about and
exact solutions of the Einstein equation with torsion. Recall that
according to Cartan \cite{Cartan25}, the torsion of a metric connection 
$\nabla$ on $(M,g)$ is the  sum of elements from
 $\Lambda^3(M), \ TM$, and a $\dfrac{n(n^2-4)}{3}$-dimensional
representation space. In dimension four, these components are called in the
physical literature the  axial vector (since $\Lambda^3(M^4)\cong TM^4$), 
the vector, and the tensor part of torsion. As will be discussed in 
Section \ref{dim4}, requiring the torsion to be skew symmetric
(i.\,e., only `axial' torsion) is rather restrictive in dimension $4$,
in particular if one imposes further mathematical conditions like parallel
torsion. Thus, most models of general relativity with torsion allow a priori
all three possibilities. Only few exact solutions to the Einstein equations
with torsion appear in the literature; these are mainly of two types,
\begin{enumerate}
\item Generalizations of classical solutions:
On the  Schwarzschild solution $(\R\x S^3, g_S)$,
one can construct a metric 
$\nabla$-Ricci flat  connection with torsion, but it is of mixed 
torsion type \cite{OMBH97},
see also Example \ref{schwarzschild}; with a rotationally symmetric Ansatz
for metric and torsion, one obtains solutions of Schwarzschild-De Sitter
\cite{B81} or Kerr-type \cite{MBG87}, again of mixed torsion.

\item Conformal changes of the flat Minkowski metric $(\R^4, g_M)$:
for example, the following Ansatz for  metric and torsion
$$
\tilde g=e^{2\omega(t)}(-dt^2+dx^2+dy^2+dz^2),\quad
H=f(t)\, dx\wedge  dy\wedge dz
$$
can be adjusted in such a way to yield an exact solution of both field
equations, this time with pure axial torsion  \cite{Lenzen84}.
 
\end{enumerate}

\noindent As far as we know, no general investigation of Einstein manifolds
with torsion was carried out. In practice,
torsion turned out to be hard to detect experimentally, since
all  tests of general relativity are based on experiments in empty space.  
Many concepts that Einstein-Cartan theory inspired are still of relevance  (see 
\cite{Hehl&M&M&N95} for a generalization with additional currents and shear,
\cite{Trautman99} for optical aspects, \cite{Ruggiero&T03} for the link to 
the classical theory of defects in elastic media). In cosmology,
Einstein-Cartan theory is again being considered in recent times
(see for example \cite{Poplawski11}).

%%%%%%%%%%%%%%%%%%%%%%%%%%%%%%%%%%%%%%%%%%%%%%%%%%%%%%%%%%%%%%%%%%%%%
\subsection{Notations and review of curvature relations}\label{subsec: notations}
%%%%%%%%%%%%%%%%%%%%%%%%%%%%%%%%%%%%%%%%%%%%%%%%%%%%%%%%%%%%%%%%%%%%%%
We end this Section by recalling a few standard identities.
Let $(M,g)$ be a Riemannian manifold and $H\in\Lambda^3(M)$.
The  metric connection with skew torsion $H$ is
$$
\nabla_X Y = \nabla^g_X Y +\frac{1}{2}H(X,Y,-).
$$
Quantities refering to the Levi-Civita connection will carry
an upper index $g$, while quantities associated with the new connection
will have an upper index $\nabla$. For example, $s^g$ and $s^\nabla$
will be the Riemannian and the $\nabla$-scalar curvatures, respectively.
>From the $3$-form $H$, we can define an associated algebraic $4$-form
$\sigma_H$, quadratic in $H$, given by
$$
2  \sigma_H\  =\  \sum_{i=1}^n (e_i \lrcorner H) \wedge (e_i \lrcorner H)
$$ 
where $e_1, e_2,\dots, e_n$ denotes an orthonormal frame of $TM$. 
The following well-known curvature identities are crucial for the
topic of this paper; they can for example
be found in \cite{FriedIvan02},  \cite{Agri06}.
We introduce the tensor
\begin{equation}\label{diff-tensor-T}
S(X,Y) :=  \sum_{i=1}^n g(H(e_i,X),H(e_i,Y)) 
= \sum_{i,j=1}^n H(e_i,X,e_j)H(e_i,Y,e_j)
\end{equation}
that measures the (symmetric part of the) difference between
the Riemannian and the $\nabla$-curvature. We normalize the
length of a 3-form $H$ as 
$\Vert H\Vert^2 = \frac{1}{6} \sum_{ij}g(H(e_i,e_j),H(e_i,e_j)).$
\begin{theorem}\label{thm-curvatures}
%---------------------------------------------------------------
%
The Riemannian curvature quantities and the $\nabla$-curvature
quantities are related by
$$\begin{array}{rcl}
R^\nabla(X,Y,Z,W) &  = &  R^g(X,Y,Z,W)+ 
\frac{1}{4} g(H(X,Y),H(Z,W)) + \frac{1}{4}\sigma_H(X,Y,Z,W)\\ & &  + \frac{1}{2}\nabla_X H (Y,Z,W) - \frac{1}{2}\nabla_Y H(X,Z,W) \medskip\\
\mathrm{Ric}^\nabla(X,Y) & = & \Ric^g(X,Y)-\frac{1}{4} 
 S(X,Y) - \frac{1}{2} \delta H (X,Y) \medskip\\
s^\nabla &  = &  s^g - \frac{3}{2} \Vert H \Vert^2
\end{array}$$
\end{theorem}
Observe that the second identity can be interpreted as the splitting
of $\Ric^\nabla$ in its symmetric and antisymmetric part. 
Where convenient, we shall use the notations
$$
S(\Ric^\nabla)\ :=\ \Ric^g(X,Y)-\frac{1}{4}  S(X,Y),\quad
A(\Ric^\nabla)\ :=\ - \frac{1}{2} \delta H (X,Y)
$$
for the symmetric and the antisymmetric part of the Ricci tensor,
respectively.
%
%---------------------------------------------------------------------------
\section{Einstein metrics with skew torsion}
%---------------------------------------------------------------------------

%---------------------------------------------------------------------------
\subsection{The variational principle}
%---------------------------------------------------------------------------
%
The standard Einstein equations of Riemannian geometry can be obtained 
by a variational argument. They are the critical points of the Hilbert 
functional
$$
g \longmapsto \int_M \left[s^g - 2 \Lambda\right]\, \mathrm{dvol}_g,
$$ 
where $\Lambda$ is a cosmological constant. Thus,
one way of obtaining Einstein equations with skew torsion is to
look for the critical points (with respect to the metric) of the 
following functional
$$
(g,H) \longmapsto \int_M \left[s^\nabla - 2\Lambda\right]\, 
\mathrm{dvol}_g = \int_M \left[s^g - \frac{3}{2} \Vert H \Vert^2  
- 2\Lambda\right]\, \mathrm{dvol}_g.
$$
For this, we will study the variation of $\Vert H \Vert^2_g$ with respect
to $g$; the torsion $H$ does \emph{not} yet need to be $\nabla$-parallel.
\begin{theorem}
%---------------
The critical points of the functional
$$
\mathcal{L}(g,H) = \int_M \left[s^\nabla - 2\Lambda\right]\, \mathrm{dvol}_g 
$$
are given by pairs $(g,H)$ satisfying the equation 
$$-S(\Ric^\nabla) + \frac{1}{2} s^\nabla g - \Lambda g= 0.
$$ 
\end{theorem}
\proof
%---------
We use the summation convention throughout the proof in order 
to increase readability.
Set $g(t) = g+t h$ and $\{e_i(t)\}$ an orthonormal basis for $g(t)$ 
such that $e_i(0) = e_i$. We have, according to our normalization,
the following identity
$$
\Vert H \Vert^2_g(t) = \frac{1}{6} g(t)(H(e_i(t),e_j(t)), H(e_i(t), e_j(t))).
$$
Taking the derivative with respect  to $t$ 
(denoted henceforth by $\del_t$) and setting $t=0$, we get
$$
\begin{array}{lcl} \del_t \Vert H \Vert^2_{g(t)}|_{t=0} & 
= & \ds \frac{1}{6} h(H(e_i,e_j), H(e_i,e_j)) 
+ \frac{1}{3} g(\del_t(H(e_i(t),e_j(t))|_{t=0}, H(e_i,e_j)) \medskip\\
   & = & \ds \frac{1}{6} h(g(H(e_i,e_j),e_k) e_k, g(H(e_i,e_j),e_l)e_l) + \\
   & & + \ds \frac{1}{3} g( H(\del_t e_i|_{t=0}, e_j) 
+ H(e_i, \del_t e_j|_{t=0}, H(e_i,e_j))\medskip\\
   & = & \ds\frac{1}{6} g(H(e_i,e_j), e_k)g(H(e_i,e_j),e_l) h(e_k,e_l) + \\
   & & + \ds\frac{1}{3} g( H(\del_t e_i|_{t=0}, e_j), H(e_i,e_j)) + \frac{1}{3} 
g( H(e_i, \del_t e_j|_{t=0}, H(e_i,e_j)). \medskip\\
\end{array}
$$
Using now the new tensor field $S$ defined by equation (\ref{diff-tensor-T}), 
we have
$$
\begin{array}{lcl}
    \del_t \Vert H \Vert^2_{g(t)}|_{t=0}    & 
= & \ds \frac{1}{6} H(e_i,e_j,e_k)H(e_i,e_j,e_l)h(e_k,e_l) 
+ \frac{1}{3} S(\del_t e_i|_{t=0}, e_i) 
+ \frac{1}{3} S(\del_t e_j|_{t=0}, e_j)\medskip\\
   & = & \ds \frac{1}{6} (S,h)_g 
+ \frac{2}{3}  S(\del_t e_i|_{t=0}, e_i) \medskip
   \ = \ \frac{1}{6} (S,h)_g 
+ \frac{2}{3}  S(g(\del_t e_i|_{t=0},e_k)e_k, e_i) \medskip\\
   & = & \frac{1}{6} (S,h)_g 
+ \ds\frac{2}{3} S(e_k, e_i) g(\del_t e_i|_{t=0},e_k)
  \end{array}
$$
If we differentiate the equality $g(t)(e_i(t), e_j(t)) =
\delta_{ij}$ with respect to $t$ and replace $t=0$, we obtain the
following equation
$$
g(\del_t e_i|_{t=0}, e_j)+g(e_i, \del_t e_j |_{t=0}) + h(e_i, e_j) = 0.
$$
Using the above identity and the fact that $T^H_g$ is symmetric, 
we then get that
$$
\begin{array}{lcl}
 \del_t \Vert H \Vert^2_{g(t)}|_{t=0} & = & \frac{1}{6} (S,h)_g 
+ \frac{1}{3} T^H_g(e_k, e_i) (g(\del_t e_i|_{t=0},e_k)+
g(e_i, \del_t e_k|_{t=0})) \medskip\\
& = & \frac{1}{6} (S,h)_g - \frac{1}{3} S(e_k, e_i) h(e_k,e_i) \medskip\\
& = & \frac{1}{6} (S,h)_g - \frac{1}{3} (S,h)_g
\ = \ -\frac{1}{6} (S,h)_g.
\end{array}
$$
Finally we can conclude that $\Vert H \Vert^2_{g(t)}$ has a first 
order Taylor expansion as follows
$$
\Vert H \Vert^2_{g(t)} = \Vert H \Vert^2_{g} - \frac{1}{6} (S, h)_g\, t 
+ o(t^2).
$$
We can then calculate the stationary points for the functional
$\mathcal{L}(g,H)$ by finding the solutions to the equation
$$
\del_t \mathcal{L}(g+th, H) = 0.
$$
Moving the derivative under the integral sign, this is then equivalent to
$$
\int_M \del_t \left[ \left(s^{g(t)}-\frac{3}{2}\Vert H \Vert^2_{g(t)} 
- 2\Lambda \right)\mathrm{dvol}_{g(t)}\right]_{|_{t=0}} =0.
$$
Recall from the classical theory that the following identities hold
$$
s^g(t) = s^g + ( \mathrm{div}(X) - (Ric^g, h)_g ) t + o (t^2),
$$ 
where $X$ is a vector field whose particular form is not important for us, and
$$
\mathrm{dvol}_{g(t)} = \mathrm{dvol}_g + \left(\frac{1}{2}(h,g)_g\, 
\mathrm{dvol}_g\right)\, t + o(t^2).
$$
Then the stationary points of our functional $\mathcal{L}$ are 
given by the equation
$$
\int_M \left[\mathrm{div}(X) - (\mathrm{Ric}^g, h)_g 
- \frac{3}{2} (-\frac{1}{6}S,h)_g  + \frac{1}{2}(s^g 
- \frac{3}{2} \Vert H \Vert^2-2\Lambda )(g,h)_g\right]\,\mathrm{dvol}_g = 0.
$$
The divergence term integrates to zero and simplifying the expression we get
$$
\int_M \left[ \, \left(-\Ric^g + \frac{1}{4}S 
+ \frac{1}{2} s^\nabla g - \Lambda g , h\right)_g \, \right] 
\mathrm{dvol}_g = 0.
$$
Noticing that $(\phantom{a},\phantom{a})_g$ is a scalar product on
the space of all symmetric 2-tensors and that $h$ is arbitrary we
can then conclude that
$$
-\Ric^g + \frac{1}{4}S + \frac{1}{2} s^\nabla g- \Lambda g = 0
$$
which is then equivalent to having $
-S(\Ric^\nabla) + \frac{1}{2} s^\nabla g - \Lambda g= 0. $

\qed

As in the Riemannian case, taking the trace of this
equality yields that $s^\nabla /2 -\Lambda=s^\nabla/n$. Consequently,
we define:
\begin{definition}
%-------------------
A triple $(M,g,H)$ is said to be `Einstein with skew torsion' or just
`$\nabla$-Einstein' if the connection $\nabla$ with torsion $H$ 
satisfies the Einstein equation
$$
S(\Ric^\nabla) = \frac{s^\nabla}{n} g.
$$
It will be called `Einstein with parallel skew torsion' if in addition
it satisfies $\nabla H=0$.
\end{definition}
In particular, $\Ric^\nabla$ does not have  to be symmetric in general. 
However, an Einstein structure  with 
\emph{parallel} skew torsion satisfies $\delta H=0$, so the symmetrization is
unnesserary in that case.

\begin{example}
All manifolds admitting a flat metric connection with skew torsion
will be trivially $\nabla$-Einstein with skew torsion. These were studied
by \'{E}. Cartan and J. Schouten \cite{CartScho26ii} who argued
(with a wrong proof) that, up
to universal cover, such manifolds are products of  a Lie group (in
the case where the torsion is parallel) or otherwise the $7$-sphere
(where the torsion is closed).  A modern,
classification-free proof using holonomy theory can be found in 
\cite{AgriFrie10}.

\end{example}

%---------------------------------------------------------------------------
\subsection{Curvature properties}
%---------------------------------------------------------------------------

One of the most important features of Riemannian Einstein spaces
is the fact that their scalar curvature is constant, and the many consequences
that follow from it.
We begin with an easy example illustrating that $\nabla$-Einstein 
manifolds with non-constant scalar curvature exist:

\begin{example}[An example on $S^3$]\label{exa-sphere}
%------------------------------------------------------
Consider the $3$-dimensional sphere $S^3$ and take $g$ to be the round
metric. Then it is well known that $(S^3,g)$ is Einstein with $s^g =
6$ and a parallelizable manifold. Let $f: S^3 \lra \R$ be any non-constant
smooth function (like the height function) and consider the
three-form $H$ given by
$H = 2 f e^1 \wedge e^2 \wedge e^3.$ 
Then the connection defined by
$\nabla = \nabla^g + \frac{1}{2}H $ is Einstein with skew torsion
and 
$$
s^\nabla = s^g - \frac{3}{2}\Vert H \Vert^2 = 6 - 6 f(x)^2,
$$
which is clearly not a constant; the example also shows that
an Einstein manifold with skew torsion can have scalar curvature
of any sign, even in the  compact case.
\end{example}
Thus, extra conditions are required to conclude the constancy of $s^\nabla$. 
We shall argue that a sufficient --- and very natural --- condition to impose
is that $\nabla H=0$, i.\,e.~the torsion of the connection is parallel.
As in the Riemannian situation, the second Bianchi identity
is the key ingredient.  The general form of the second Bianchi identity
of anly linear connection may be found in the standard reference
\cite{KobaNomi69}; however, we need it in our special situation, 
where it is easier to derive it directly than to get it by specializing
the general formula.
We denote by the symbol $\stackrel{XYZ}{\sigma}$
the cyclic sum over $X,Y$, and $Z$.

\begin{proposition}[Second Bianchi identity]
%-------------------------------------------
Let $(M,g,H)$ be a Riemannian manifold equipped with a connection
$\nabla$ with skew torsion $H$ such that $\nabla H = 0$. Then
$$\stackrel{WXY}{\sigma} \nabla_W R^\nabla(X,Y)Z = \stackrel{WXY}{\sigma} \left( R^g(W, H(X,Y))Z - \frac{1}{2}R^g(W,X)H(Y,Z)) \right) $$
\end{proposition}

\proof 
%------
For the $(1,3)$-curvature tensor quantities the following
identity holds
$$R^\nabla(X,Y)Z = R^g(X,Y)Z + \frac{1}{2}H(H(X,Y),Z)+ \frac{1}{4}H(H(Y,Z),X)+\frac{1}{4}H(H(Z,X),Y)$$
Set $$R^H(X,Y)Z = \frac{1}{2}H(H(X,Y),Z)+
\frac{1}{4}H(H(Y,Z),X)+\frac{1}{4}H(H(Z,X),Y)$$ and notice that
since $\nabla H = 0$ then also $\nabla R^H = 0$. Then it is easy to
check that
$$\nabla_W R^\nabla (X,Y)Z = \nabla^g R^g(X,Y)Z - \frac{1}{2}(R^g(H(W,X),Y)Z + R^g(X,H(W,Y))Z + R^g(X,Y)H(W,Z))$$
Now we just need to take the cyclic permutation, use the second
Bianchi identity for $R^g$ and the proposition follows.
\qed

\begin{corollary}
%----------------
Assume $\nabla H=0$.
The divergence of the $\nabla$-Ricci tensor is proportional to the
derivative of the $\nabla$-scalar curvature, more precisely:
$$\delta \Ric^\nabla +\frac{1}{2} d{s^\nabla} =0.$$
\end{corollary}

\noindent\emph{1st proof} ---
%-----------------------------
Taking traces in the second Bianchi identity for $R^\nabla$
we immediately get the following equation

$$\begin{array}{l}-\nabla_X R^\nabla(e_i,e_j,e_j,e_i) + \nabla_{e_i} R^\nabla(e_j,X,e_i,e_j) + \nabla_{e_j} R^\nabla (e_i,X, e_j,e_i) = \medskip\\
\phantom{-\nabla_X R^\nabla} = R^g(X,H(e_i,e_j), e_i, e_j) +
R^g(e_j, H(X,e_i),e_i,e_j) + R^g(e_i, H(e_j,X),e_i,e_j)
\medskip \\ \phantom{-\nabla_X R^\nabla +}
-\frac{1}{2}R^g(X,e_i,H(e_j,e_i),e_j)-\frac{1}{2} R^g(e_i,e_j,
H(X,e_i),e_j) \end{array}$$ which then simplifies to
$$\begin{array}{l}-d s^\nabla (X) - 2 \delta^\nabla \Ric^\nabla(X)  = \\ \phantom{-d s^\nabla } = R^g(X, H(e_i, e_j),e_i,e_j)+
\frac{1}{2}R^g(X,e_i,H(e_i,e_j),e_j) +\frac{5}{2}\Ric^g(H(X,e_i),e_i)\medskip\\
\phantom{-d s^\nabla }  =  H(e_i,e_j,e_k)R^g(X, e_k,e_i,e_j)+
\frac{1}{2}H(e_i,e_j,
e_k)R^g(X,e_i,e_k,e_j)+\frac{5}{2}H(X,e_i,e_k)\Ric^g(e_k,e_i)\end{array}$$
Now the first two terms on the right-hand-side simplify because $H$
is antisymmetric and $R^g$ satisfies the first Bianchi identity and
the last term vanishes because $H$ is antisymmetric and $\Ric^g$ is
symmetric. Finally, observe that since $H$ is totally antisymmetric, the
$\nabla$-divergence of any symmetric $(0,2)$-tensor is the same as the
usual divergence (\cite{Agri06}, Proposition A.2), hence
$\delta^\nabla \Ric^\nabla = \delta \Ric^\nabla$.

\qed

\medskip
\noindent\emph{2nd proof} ---
%-----------------------------
Let $S$ be the symmetric 2-tensor introduced before
that satisfies  $\Ric^\nabla = \Ric^g -S/4$. 
Notice that $S$ is $\nabla$ parallel since it is a composition 
of two parallel tensors. Then
$$
\delta^\nabla \Ric^\nabla = \delta^\nabla \Ric^g,
$$
and this is again $\delta \Ric^g$ by the preceding comment on divergences.
But from the Riemannian case, we know that $\delta^g \Ric^g 
= -\frac{1}{2}d s^g$. Since $H$ is parallel, $d(\Vert H \Vert^2) = 0$, 
hence the relation between the $\nabla$- and the Riemannian scalar curvature
implies $ds^g = ds^\nabla$.

\qed

Using the corollary of the second Bianchi identity we can prove the
following.

\begin{proposition}\label{Prop: scal-constant}
%----------------------------------------------
Assume $\nabla H=0$. If $\delta\Ric^g=0$, the
scalar curvatures $s^\nabla$ and $s^g$ are constant. In particular,
this holds if $(M,g,H)$ is Einstein with parallel skew torsion $H$.
\end{proposition}
\proof 
%--------
The first claim is immediate. For the $\nabla$-Einstein case,
taking the divergence on both sides of the equation
$\Ric^\nabla = \frac{s^\nabla}{n}g $ and using the proposition above
we get that $-\frac{1}{2}d s^\nabla = \frac{1}{n} d s^\nabla$.
Therefore $d s^\nabla$ = 0 and $s^\nabla$ is constant.
But as observed before, a parallel torsion form has constant length, hence
the relation between the scalar curvatures (see Theorem 
\ref{thm-curvatures}) implies that the Riemannian scalar curvature will be
constant as well.

\qed
\begin{remark}
%--------------
There are other situations with $\nabla H=0$ where one can conclude that
the scalar curvatures $s^\nabla$ has to be constant: For example,
this holds if $M$ is spin and if there exists a non trivial parallel 
spinor field $\psi$, $\nabla\psi=0$ (see \cite[Cor. 3.2]{FriedIvan02}). 
Since any $G_2$ manifolds with a characteristic connection $\nabla$ admits a
$\nabla$-parallel spinor (see Section \ref{section:G2}) and many 
of them are known
to have parallel torsion, many examples of this kind that are not
Einstein with skew torsion exist.
\end{remark}
\begin{corollary}
%----------------
Any complete connected Riemannian manifold $(M,g,H)$ that is
Einstein with parallel skew torsion $H$ and with positive
scalar curvature $s^\nabla>0$ is compact and has finite first
fundamental group $\pi_1(M)$. 
\end{corollary}
\proof
The crucial point is that the assumption of the
Bonnet-Myers Theorem has to hold, i.\,e.~the inequality
$\Ric^g(X,X)\geq c \|X\|^2$ for some positive constant $c$
and all $X\in TM$. But this is easy: 
$$
\Ric^g(X,X)\ = \ \Ric^\nabla(X,X)+\frac{1}{4}S(X,X)\ =\
\frac{s^\nabla}{n}\|X\|^2 + \frac{1}{4}S(X,X)\ \geq \frac{s^\nabla}{n}\|X\|^2,
$$
since $S$ is a non-negative tensor by definition. All claims now follow
from the classical Riemannian results. Observe that is is not
necessary to specify further with respect to which connection
completeness is meant: A metric connection with skew torsion on a
Riemannian manifold $(M,g)$ is complete if and only if the Levi-Civita
connection is complete, because their geodesics coincide.

\qed

We describe now how our notion of Einstein with parallel skew
torsion is consistent with the algebraic decomposition of the
curvature tensor.
Let $\mathcal{C}M$ be the space of all symmetric curvature tensors
on $TM$. Consider the Bianchi map
$$
b:  \mathcal{C}M  \lra  \mathcal{C}M,\quad
    R_{abcd}  \lms  R_{abcd}+R_{bcad}+R_{cabd}.
$$
It is well known (see \cite{Bess87i} for details) that
$S^2 (\Lambda^2 M) = \ker b \oplus \im\, b$ and that $\im\, b = \Lambda^4
M$.
Suppose now that our Riemannian manifold with torsion $(M,g, H)$ is
such that $H$ is $\nabla$ parallel. Then the Riemannian curvature
tensor simplifies to
$$R^\nabla(X,Y,Z,W) = R^g(X, Y, Z, W) + \frac{1}{4}g(H(X,Y),H(Z,W)) +\frac{1}{4} \sigma_H(X,Y,Z,W).$$
Observe that in this case $R^\nabla$ is indeed in $S^2(\Lambda^2
M)$. We have the following proposition:

\begin{proposition}
%------------------
Let $(M,g,H)$ be such that $\nabla H = 0$. Then $R^\nabla$
decomposes under the Bianchi map as $R^\nabla = R^\nabla_k +
R^\nabla_i$, where $R^\nabla_k$ lies in the kernel of the Bianchi map and $R^\nabla_i$ in its image, with
$$R^\nabla_k(X,Y,Z,W) = R^g(X,Y,Z,W)+\frac{1}{4}g(H(X,Y),H(Z,W)) -\frac{1}{12}
\sigma_H(X,Y,Z,W)$$ and $$R^\nabla_i(X,Y,Z,W) =
\frac{1}{6}\sigma_H(X,Y,Z,W)$$
\end{proposition}

\proof Notice that $\stackrel{XYZ}{\sigma}{g(H(X,Y),H(Z,W)} =
\sigma_H(X,Y,Z,W)$ and $\stackrel{XYZ}{\sigma}\sigma_H(X,Y,Z,W) =
3\sigma_H(X,Y,Z,W)$. Since $R^g$ is  in $\ker b$ then so is
$R^\nabla_k$. Notice also that $\sigma_H \in \Lambda^4 M$.

\qed

Following the classical theory we can then decompose the
$\nabla$-curvature tensor as
$$
R^\nabla = W^\nabla + \frac{1}{n-2} (Z^\nabla \boxtimes g) +
\frac{s^\nabla}{n(n-1)} g \boxtimes g + \frac{1}{6} \sigma_H
$$ 
where $\boxtimes$ denotes the Kulkarni-Nomizu product, 
$Z^\nabla =\Ric^\nabla - \frac{s^\nabla}{n}g $ is the
trace-free part of the $\nabla$-Ricci tensor and $W^\nabla$, which
we shall call the $\nabla$-Weyl tensor, can be explicitly written as
$$
\begin{array}{lcl} W^\nabla(X,Y,Z,W) & = & W^g(X,Y,Z,W) + \dfrac{1}{4}g(H(X,Y), H(Z,W)) + \sigma_H(X,Y,Z,W)\\
& & + \dfrac{1}{4(n-2)} ( g(H(e_i,X),H(e_i,W))g(Y,Z) +
g(H(e_i,Y),H(e_i,Z))g(X,W) \\
& & \phantom{+ \dfrac{1}{4(n-2)}}- g(H(e_i,X),H(e_i,Z))g(Y,W) -
g(H(e_i,Y),H(e_i,W))g(X,Z))\\
& & - \dfrac{3\Vert H \Vert^2}{2(n-1)(n-2)}
(g(X,W)g(Y,Z)-g(X,Z)g(Y,W))
\end{array}
$$
Note that $W^\nabla$ is traceless. This can be checked by direct
calculation but, of course, it also follows from the general theory.
We conclude that our previous definition of $(M,g,H)$ being
Einstein with skew torsion is equivalent to $Z^\nabla =0$, as it should;
let us emphasize that this relies again strongly on the property $\nabla H=0$.

Next, we want to clarify the relation between $\nabla$-Einstein
and Riemannian Einstein manifolds. In a given dimension,
the algebraic form of the 
difference tensor $S(X,Y)$ of the curvatures as defined in 
equation \ref{diff-tensor-T}
decides whether a Riemannian Einstein metric will yield a skew Einstein
structure or not.
\begin{definition}
%------------------
On a Riemannian manifold $(M,g)$, a $3$-form $H$ will be called `of 
Einstein type' if the difference tensor $S(X,Y) := \sum_i g(H(e_i,X),
H(e_i,Y))$ is proportional to the metric $g$.
\end{definition}
\begin{proposition}\label{Prop:relation}
%----------------------------------------
Let $(M,g)$ be a Riemannian manifold, $H$ a $3$-form 
written in a local orthonormal frame $e^1,\dots,e^n$ of $T^*M$,
$
H = \sum_{ijk} H_{ijk}\, e^i \wedge e^j \wedge e^k.
$ 
Then  $H$ is of Einstein type if it satisfies the following conditions:
\begin{enumerate}
\item no term of the form $H_{ija} e^i\wedge e^j \wedge e^a + H_{ijb} e^i\wedge e^j \wedge
e^b$ with $a\neq b$ occurs;
\item if $i$ and $j$ are two indices in $\{1,\dots, n\}$ then the
number of occurrences of $i$ and $j$ in $H$ coincides;
\item if $\{i,j,k\}$ and $\{a,b,c\}$ are two sets of indices then $H_{ijk}^2 =
H_{abc}^2$.
\end{enumerate}
\end{proposition}

\proof It can be checked by direct calculation
that, when writing $S$ in matrix form for the frame $\{e_i \}$,
condition (1) guarantees $S$ to be diagonal, while conditions (2) and (3)
ensure that $S$ is indeed a multiple of the identity matrix.

\qed

\begin{remark} 
%--------------
Proposition $\ref{Prop:relation}$ yields an easy procedure for
producing further examples of $\nabla$-Einstein metrics with non
constant scalar curvature (beyond the one given in Example \ref{exa-sphere})
for all manifolds that are parallelizable and carry an Einstein metric 
(for example $S^7$ or compact semi-simple Lie groups).
\end{remark}

\begin{figure}
\begin{center}
\psfrag{I}{$e_{123}$}
\psfrag{II}{$e_3(e_{15}+e_{24})$}
\psfrag{III}{$e_{126}+e_{135}+e_{234}$}
\psfrag{IV}{\!\!\!${}^{*}\ e_{123}+e_{456}$}
\psfrag{V}{\!\!\!\!\!${}^{*}\ e_{135}+e_{146}+e_{256}+e_{234}$}
\psfrag{6}{VI. $e_1(e_{23}+e_{45})+e_{267}$}
\psfrag{7}{VII. $e_{123}+e_{456}+e_7(e_2+e_5)(e_3+e_6)$}
\psfrag{8}{VIII. $e_1(e_{23}+e_{45}+e_{67})$}
\psfrag{9}{IX.$^{**}$ VIII + $e_{246}$,}
\psfrag{10}{X. VII + $e_{147}$}
\psfrag{11}{XI. VIII + $e_2(e_{46}- e_{57})$}
\psfrag{12}{XII. $e_1(e_{45}+e_{67})+e_2(e_{46}-e_{57})+ e_3(e_{47}+e_{56})$,}
\psfrag{13}{XIII.$^{**}$ XII $- e_{123}$}
\includegraphics[width=8cm]{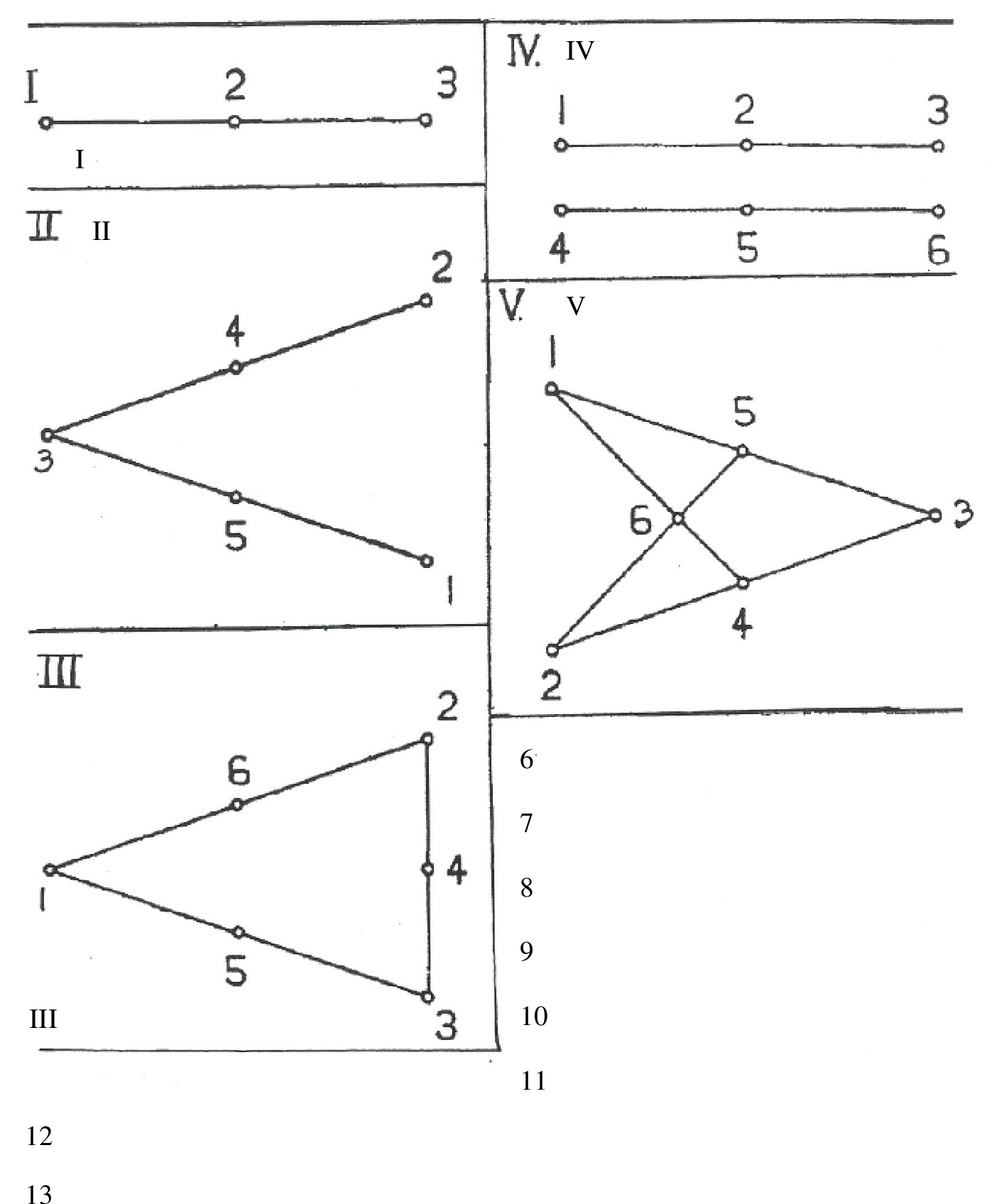}
\end{center}
\caption{Normal forms of $3$-forms up to dimension $7$ (see Example 
\ref{exa:normal-forms}); the picture is
taken from the original article of J.~Schouten from 1931 \cite{Schouten31}
while for $n=7$ we follow \cite{Westwick81}. Any three numbers $a,b,c$ linked
by lines represent a summand $e_a\wedge e_b\wedge e_c$.}
\label{normal-forms}
\end{figure}

\begin{example}\label{exa:normal-forms}
%--------------------------------------
In 1931, J.~Schouten described the normal forms of 
$3$-forms up to dimension $7$ \cite{Schouten31}, i.\,e.~representatives of the 
$\mathrm{GL}(n,\R)$-orbits for $n\leq 7$, see Figure \ref{normal-forms} 
(the complex classification is different; see also \cite{Westwick81} 
for a modern account of the real classification). 
The normal forms marked with one resp.~two stars are the representatives of
generic $3$-forms, i.\,e.~the ones with
dense $\mathrm{GL}(n,\R)$ orbit for $n=6$ resp.~$n=7$ (see for example
\cite{Vanzura&W12}).
One checks
by a direct computation that only the following $3$-forms 
are of Einstein type in the given dimensions:
\begin{enumerate}
\item Type I. in dimension $3$,
\item Types IV. and V. in dimension $6$,
\item Type XIII. in dimension $7$.
\end{enumerate}

In particular, Riemannian Einstein manifolds $(M,g)$ will never be
$\nabla$-Einstein in dimensions $4$ and $5$. We will  see later 
concrete examples where Proposition $\ref{Prop:relation}$ can be 
applied to construct $\nabla$-Einstein metrics. 
\end{example}

\begin{example}\label{schwarzschild}
%------------------------------------
Proposition $\ref{Prop:relation}$ shows that there does not
exist a three-form that will make
the Schwarzschild metric $\nabla$-Einstein; in particular, the connection
constructed in \cite{OMBH97} has to be of mixed torsion type.
\end{example}
%
%-----------------------------------------------------------------------
\subsection{Examples and construction of $\nabla$-Einstein manifolds}
%-----------------------------------------------------------------------

\subsubsection{The case $n=4$}\label{dim4}
%----------------------------------------
In dimension $4$, an alternative notion of Einstein with skew torsion
was investigated by the second author in \cite{Ferr10, Ferr11},
based on the phenomenon of self-duality. The idea was to
consider the decomposition of the curvature operator
$\mathcal{R}: \Lambda^2 \rightarrow \Lambda^2$ in terms of the splitting
of the bundle $\Lambda^2 = \Lambda_+\oplus \Lambda_-$ into self-dual
and antiself-dual parts. By making the analogy with the standard 
definition of Einstein metric, one can set the upper-right block 
in this decomposition to vanish, that is, the following definition was taken.
In order not to confuse it with the Einstein definition used in this
article, we will call it \emph{duality-Einstein with skew torsion}, 
in the sense that it is the Einstein characterisation based on 
(self-)duality.
\begin{definition}
%-------------------
A triple $(M^4,g,H)$ is said to be duality-Einstein with skew torsion if
$$
Z^\nabla + S \left( \nabla\! *\! H \right)+ \dfrac{*dH}{4}  g = 0,
$$
where $S$ denotes the symmetrization of a tensor 
and $Z^\nabla$ is the symmetric trace-free part of $\Ric^\nabla$.
\end{definition}
This definition depends a priori on the choice of orientation, but
it can be proved that for a compact manifold, this choice is irrelevant,
since the metric dual of $*H$ is a Killing vector field. Among other
results, it was shown that the Hitchin-Thorpe inequality 
$2 \chi \geq 3 |\tau|$ holds again, giving a severe
topological obstruction.
Since, in four dimensions, $*H$ is a $1$-form, it is not difficult
to establish that the notion of duality-Einstein manifold with skew torsion
is equivalent to that of Einstein-Weyl manifold, which has been a
subject of intense study in the past \cite{CaldPede99}. 
Our new definition of Einstein manifolds with parallel
skew torsion breaks away from Einstein-Weyl geometry for $n\neq 4$.
This can be seen for example from the scalar curvature: 
while it has to be constant in our situation by 
Proposition \ref{Prop: scal-constant}, it has only
constant sign in Einstein-Weyl geometry \cite{PedeSwan93i}.

In general, the notions of $\nabla$-Einstein and duality-Einstein
differ. But with the additional assumption  that the torsion is parallel,
they coincide: first, observe that in $n=4$, condition
$\nabla H=0$ is equivalent to $dH=0$ by 
Proposition \ref{Prop: Bianchi} and the fact that  $\sigma_H=0$ in
dimension four for purely algebraic reasons. Also,
$\nabla *H=0$ as well, so the  duality-Einstein equations reduces 
to $Z^\nabla=0$, and this is  our definition of an Einstein structure 
with parallel skew torsion.
In fact, more holds: since $*H$ is a parallel vector field,
such a manifold will have to be $\nabla$-Ricci-flat. 
A result similar to the following (but under different assumptions)
was proved in the second author's dissertation \cite{Ferr10}.

%-----------------------------------------------------------------
\begin{theorem}\label{Theo: Bismut-Einstein-Skew}
%-----------------------------------------------------------------
%
If $(M,g,H)$ is a $4$-dimensional complete, simply connected
manifold that it is  Einstein with parallel skew torsion $H\neq 0$, 
then it is isometric to $\mathbb{R}\times S^3$ with the standard product
metric.
\end{theorem}

\proof 
%--------
Define $V=*H$, which is of constant length, $\nabla$-parallel, and non-zero by
assumption. One easily sees that $\nabla^g V=0$. The vector field $V$
is complete, because it is a Killing vector field on a complete Riemannian
manifold. By the de Rham theorem, we
have a reduction of $TM$ under the holonomy group and that
$M$ is isometric to $\R\times N$, where $N$ is a 3-dimensional
manifold such that $TN$ is the orthogonal complement of $\{V\}$.
Now $\mathrm{Ric}^\nabla = 0$, and any torsion tensor in dimension $3$
is of Einstein type (see Example \ref{exa:normal-forms}), so
$N$ is Riemannian Einstein as well. Furthermore, $s^\nabla=0$ implies
 $s^g = 3/2 \|H\|^2 >0$, 
i.\,e.~$N$ is a Riemannian Einstein $3$-manifold of positive scalar 
curvature, hence the sphere. We can conclude that $M$ is 
isometric to $\R\times S^3$.

\qed

\subsubsection{Lie groups}
%-------------------------
Classical examples of manifolds where skew torsion arises naturally
are  Lie groups $G$ equipped with a bi-invariant inner product $g$
on the corresponding Lie algebra $\g$ (see \cite{KobaNomi69}, for
example).
There is a one-to-one correspondence between the set of all
bi-invariant connections and the space of all bilinear functions
$\alpha:\g\times\g \lra \g$ which are ad-invariant; in particular,
$\nabla^t_X(Y) = t [X,Y]$ defines a one-parameter family of metric
connections with parallel skew torsion $(2t-1) [X,Y]$. 
The value $t=1/2$ corresponds to the Levi-Civita connection,
while $t=0$ and $t=1$ are the flat $\pm$-connections introduced by
Schouten and Cartan. A routine calculation shows:
\begin{lemma}
%--------------------
Let $G$ be a Lie group equipped with a bi-invariant metric. Consider
the 1-parameter family of connections with skew torsion
$\nabla^t_X Y = t[X,Y].$ It has Ricci curvature
$$
\Ric^{\nabla^t}(X,Y) \ =\ 4(t-t^2)\,\Ric^g(X,Y).
$$
Hence, for $t=0,1$, the connection $\nabla^t$ is $\nabla^t$-Einstein 
(it is in fact flat), and for $t\neq 0,1$, $\nabla^t$ is 
Einstein if and only if $\nabla^g$ is Einstein.
\end{lemma}
It is well known  that for any compact semisimple Lie
group, the Killing form $K$ is negative
definite and therefore $-K$ defines a bi-invariant metric. This
metric is clearly always Einstein in the usual sense.
However,  not all Lie groups admit Einstein metrics;
this happens for example for $G = U(1)\times SU(2)$.
If $(G,g)$ was Einstein, the Euler characteristic of $G$ would
be given by \cite{Bess87i}
$$
\chi(G) =\dfrac{1}{8 \pi^2} \int_G \left( \dfrac{s^2}{24} 
+ \Vert W \Vert^2 \right) \mathrm{dvol}_g.
$$
Since $\chi(G) =0$, this would  mean that both the scalar and
the Weyl curvature vanish and therefore $G$ would be flat with
respect to the Levi-Civita connection, which cannot happen since its
universal cover is $\mathbb{R}\times S^3$ and not $\mathbb{R}^4$.

%----------------------------------------------------------------------------
\subsubsection{Almost Hermitian manifolds}
%-----------------------------------------------------------------------------
%
We will know take a look at almost Hermitian manifolds. 
By a result of Friedrich and Ivanov \cite{FriedIvan02}, any almost 
Hermitian manifold
of Gray-Hervella type $W_1\oplus W_3\oplus W_4$ \cite{GrayHerv80} 
admits a unique Hermitian connection $\nabla$ with skew-symmetric torsion
(i.\,e.~a connection satisfying $\nabla g=\nabla J=0$, see \cite{Gaud97}
for a survey).
The connection is called the \emph{characteristic connection} of the
almost Hermitian structure. 

\begin{example}\label{CYT}
%--------------------------
By definition, a 
\emph{Calabi-Yau connection with torsion} (CYT connection for short) is 
a Hermitian $2n$-dimensional manifold suich that the
restricted holonomy of its characteristic connection lies in $\SU(n)$. 
In \cite{Grantcharov08}, the authors investigated such connections on
principal $T^{2k}$-bundles over compact K\"ahler manifolds.
By \cite{Grantcharov08}, Proposition 5 and equation (11), 
one sees that the existence of a CYT connection implies that the bundle has
to be $\Ric^\nabla$-flat (the connection of interest corresponds to $t=-1$).
By a sophisticated topological construction,
the authors construct several series of examples on 
$(k-1)(S^2\x S^4)\# k(S^3\x S^3)$ for all $k\geq 1$.
Examples on $S^3\x S^3$ were shown to exist previously in \cite{Fino04}.
\end{example}

\begin{example}
%------------------
An almost Hermitian (non-K\"ahlerian) manifold is called 
\emph{nearly K\"ahler} if it satisfies $(\nabla^g_X J)(X)=0$,
these are precisely the manifolds of type $W_1$.
Nearly K\"ahler manifolds were introduced and extensively studied by
A.~Gray \cite{Gray70, Gray76}. The characteristic
connection is then called the Gray connection. 
It is a non-trivial result that its torsion is $\nabla$-parallel, \cite{Kiri77,
AlexFrieScho05}. Furthermore, A. Gray showed that
any $6$-dimensional nearly K\"{a}hler manifold is Einstein \cite{Gray76}.
A routine computation in the theory of nearly K\"ahler manifolds
shows that $\displaystyle S(X,Y) =  \frac{4 s^g}{30} g(X,Y)$ and
hence $\Ric^{\nabla} = \dfrac{2 s^g }{15} g$:
Any nearly K\"{a}hler manifold is Einstein with parallel skew torsion
with respect to the Gray connection. This has been implicitly
known  in the community and was noted several times in the literature,
for example, in \cite{AlexFrieScho05}, \cite{Agri06}.
In \cite{Butr05}, it was shown that a
$6$-dimensional locally homogeneous nearly K\"{a}hler  manifold has to
be one of the following: $S^6$, $S^3\times S^3$ or the twistor
spaces for $S^4$ and $\mathbb{C}P^2$ ($\mathbb{C}P^3$ and the flag
manifold $F(1,2)$, resp.).
\end{example}

N.~Schoemann investigated in \cite{Scho07} almost Hermitian
$6$-manifolds admitting a  characteristic connection with parallel torsion.
He discovered that there is a class of manifolds whose curvature tensor
is the same as in the nearly K\"ahler case, and hence they will be
Einstein with skew torsion. Let us describe them in more detail.

\begin{theorem}[$\mathrm{\cite[Section\, 4.2]{Scho07}}$]
%-----------------------------------------------------
Let $(M^6,g)$ be an almost Hermitian manifold admitting a 
characteristic connection $\nabla$ with parallel torsion $H$
such that its reduced holonomy $\mathrm{Hol}_0(\nabla)$ lies inside 
$\mathrm{SO}(3)$. Then $(M^6,g)$ is of type $W_1\oplus W_3$,
locally isomorphic to an isotropy irreducible homogeneous space,
and  the curvature transformation is of the form
$R^\nabla: \ \Lambda^2(\R^6)\rightarrow \mathfrak{so}(3),\ 
R^\nabla(X,Y)=\lambda\mathrm{pr}_{\mathfrak{so}(3)}(X\wedge Y),\ \lambda\in\R$.
\end{theorem}

In particular, the nearly K\"ahler structure on $S^3\times S^3$ is of
this type, and thus we can conclude without calculation that
any  almost Hermitian manifold with the same curvature transformation
is $\nabla$-Einstein with parallel skew torsion. Examples are
given in Section 4.5 of \cite{Scho07} and in \cite{AlexFrieScho05}: 
Besides a larger class of
almost Hermitian structures on $S^3\times S^3$ that includes the
nearly K\"ahler case, examples can be constructed on $\mathrm{SL}(2,\C)$
(viewed as a real $6$-manifold, $\Ric^\nabla = -\dfrac{1}{3}\|H\|^2g$),
$\mathrm{SU}(2)\ltimes \R^3$, and a nilpotent group $N^6$.

%-------------------------------------------------------------------------
\subsubsection{Almost contact metric manifolds}
%-------------------------------------------------------------------------
%
We shall now investigate contact manifolds in greater detail:
We will be able to construct large classes of new Einstein structures
with torsion from them.

\begin{definition}
%-----------------
We say that $(M,g)$ carries an
almost contact metric structure if it admits a (1,1)-tensor field $\varphi$
and a vector field $\xi$ with dual
form $\eta$ such that $\varphi^2 = -Id + \eta \otimes \xi$
and if $g$ is $\varphi$-compatible, i.\,e.
$$g
(\varphi(X),\varphi(Y)) = g(X,Y)-\eta(X)\eta(Y).
$$ 
\end{definition}

\noindent
Almost contact manifolds have two naturally associated tensors, the
fundamental form $F$ defined by
$F(X,Y) = g(X,\varphi(Y))$ and the Nijenhuis tensor $N$ given by
$$N(X,Y) = [\varphi(X),\varphi(Y)]-\varphi([X,\varphi(Y)])-\varphi([\varphi(X),Y]) + \varphi^2([X,Y])+d\eta(X,Y)\xi.$$

We say that the almost contact metric manifold M is normal if $N=0$
and that it is contact metric if $2F = d\eta$. A contact metric
manifold that is also normal is called a Sasaki manifold.

Again,  Friedrich and Ivanov characterised the almost contact metric 
structures that admit a metric connection with skew torsion
satisfying $\nabla g = \nabla \eta = \nabla \varphi = 0$:
These are precisly the manifolds for which $N$ is a $3$-form
and  $\xi$ is a Killing field \cite{FriedIvan02}, the connection
is unique and again called the characteristic connection.
A class of manifolds satisfying these conditions are the Sasaki manifolds. 
In this particular case, the torsion tensor simplifies to
$H = \eta\wedge d\eta$ and also we have that $H$ is
$\nabla$-parallel.

\begin{lemma} \label{Prop: Sasaki-Ricci-Flat}
%--------------------------------------------
Let $(M,g,\varphi,\xi,\eta)$ be an  almost contact metric 
structure admitting a characteristic connection. If $M$ is $\nabla$-Einstein, 
then it is $\nabla$-Ricci flat.
\end{lemma}

\proof This follows immediately from the fact that $\nabla \xi = 0$
implies $\Ric^\nabla(\xi,\xi) = 0$.

\qed

\noindent
We will now investigate the so-called \emph{Tanno deformation} of a Sasaki
metric. The following results appears implicity in the articles of
S. Tanno \cite{Tann68} and of E. Kim and T. Friedrich
\cite{KimFrie00}, see also J. Becker-Bender's PhD thesis for a more
explicit formulation and proof \cite{BB12}.

\begin{theorem}[$\mathrm{\cite{Tann68, KimFrie00}}$] 
\label{Theo: Tanno-Sasaki} 
%-----------------------------------------------------
Let $(M,g, \xi,\varphi, \eta)$ be a Sasaki manifold of dimension $2k+1$ and $g_t =
t g + (t^2-t) \eta \otimes \eta$, $t>0$, a deformation of the metric
in the direction of $\xi$. Then $(M,g_t,\xi_t, \varphi, \eta_t)$,
were $\xi_t = \frac{1}{t} \xi$ and $\eta_t = t \eta$, is also a
Sasakian manifold. The Levi-Civita connection for the metric $g_t$
is given by
$$\nabla^{g_t} = \nabla^g + (1-t) (\eta\otimes \varphi + \varphi\otimes
\eta).$$ For the Ricci tensor we have the following expression
$$\Ric^{g_t} = \Ric^g + 2(1-t) g - (1-t)(2kt+2k-1)\eta\otimes\eta.$$
\end{theorem}

\noindent
We will now show that if we start with an Einstein-Sasaki or an 
$\eta$-Einstein-Sasaki manifold, i.\,e.~a Sasaki manifold whose
Riemannian Ricci tensor has the form $\Ric^g = a g + (b-a)\eta\otimes \eta$
(with one extra condition),
there is a parameter $t$ in the Tanno deformation such that the
characteristic connection is Einstein with skew torsion.
Recall that for $n\geq 5$, the coefficients $a$ and $b$ have to be
constant for an $\eta$-Einstein-Sasaki manifold \cite{Blair2010}.

\begin{theorem}\label{Theo: Einstein-Sasaki}
%---------------------------------------------
Let $(M^{2k+1},g,\xi,\varphi, \eta)$ be an Einstein-Sasaki manifold
or an $\eta$-Einstein-Sasaki manifold satisfying $2k+1> b-a$.
Then there exists a parameter  $t>0$ for which the
Tanno deformation $(M,g_t,\xi_t, \varphi,
\eta_t)$ equipped with its characteristic connection
$$
\nabla^t = \nabla^{g_t} + \frac{1}{2} \eta_t \wedge d\eta_t
$$ 
is $\nabla^t$-Ricci flat.
\end{theorem}

\proof We start by observing that if $H^t$ is the torsion of
$\nabla^t$ then $$H^t = \eta_t \wedge d\eta_t = t^2 d\eta \wedge
\eta = t^2 H.$$ Also since $M$ is Sasaki, $d\eta = 2 F$ where $F$ is
the fundamental 2-form, that is, $F(X,Y) = g(X,\varphi(Y))$. Notice
that then we can write $H$, viewed as a (1,2)-tensor, as
$H(X,Y)= 2 F(X,Y) \xi$. Consider an adapted orthonormal basis for
$g$, say, $\{a_1,\varphi(a_1),\dots, a_k,\varphi(a_k), \xi\}$. Then
the set $$\{ t^{-1/2} a_1, t^{-1/2}\varphi(a_1), \dots, t^{-1/2}a_k,
t^{-1/2}\varphi(a_k), t^{-1}\xi\}$$ is an adapted orthonormal basis
for the metric $g_t$. For notational convenience we will relabel
this basis by $\{e_i(t), i=1,\dots,2k+1\}$. We need to analyze the
tensor
$$
S^t(X,Y) \ =\  g_t (H^t(e_i(t), X), H^t(e_i(t), Y)).
$$
Expanding the expression for $g_t$ and $H^t$ we have
%\marginpar{\qquad please check formula}
$$
S^t(X,Y) = t^5 \sum_{i=1}^{2k+1} g(H(e_i(t),X), H(e_i(t),
Y)+ 4(1-t) \eta(H(e_i(t),X))\eta(H(e_i(t),Y))
$$ 
and using the fact
that $H(X,Y) = 2 g(X,\varphi(Y))\xi$, this simplifies to
$$
S^t(X,Y)  = 4 t^6 \sum_{i=1}^{2k+1} g(e_i(t),\varphi(X))g(e_i(t),\varphi(Y)). 
$$
Given the particular expression of our adapted basis, and recalling
that $\varphi(\xi) = 0$, we can yet write $S^t$ as 
$$
S^t(X,Y) = 4 t^5\sum_{i=1}^{2k} g(e_i,\varphi(X))g(e_i, \varphi(Y)).
$$ It is easy to
check that $S^t(a_i,a_i) = S^t(\varphi(a_i),\varphi(a_i)) = 4 t^5$,
for $i=1,\dots,k$ and that all other terms vanish. Then $S^t$ can be
written in matrix form as $S^t = 4 t^5 \mathrm{diag}(1,\dots,1,0)$.
Observe also that $\eta\otimes\eta$ is given in this basis by
$\mathrm{diag}(0,\dots,0,1)$. Hence, the 
Ricci tensor is given by the expression
$$
\Ric^{\nabla^t} = Ric^g + 2 (1-t)g - (1-t) (2kt+2k+1) \eta\otimes\eta 
-\frac{1}{4} S^t.
$$
If $(M,g)$ is Einstein, then
$\nabla^t$ will be Einstein if and only if the matrix
$$
-t^5 \mathrm{diag}(1,\dots,1,0) -
(1-t)(2kt+2k+1)\mathrm{diag}(0,\dots,0,1)
$$ is a multiple of the
identity. This happens if and only if
$$
f_1(t):= t^5-(1-t)(2kt+2k+1) = 0.
$$ 
By the Intermediate Value Theorem we can
conclude that such a solution exists for some $t>0$
(in fact, one checks numerically that for $k=2,\ldots,10$, the unique
positive solution lies in the intervall $[9/10,1[$). Thus $\nabla^t$
is Einstein. 

In case $(M,g)$ is only $\eta$-Einstein, one deduces similarly
that one needs to solve the equation 
$$
f_2(t):= t^5-(1-t)(2kt+2k+1) +(b-a) = 0.
$$
The assumption $2k+1>b-a$ then implies that $f(0)<0$, hence
the Intermediate Value Theorem garantees again the existence of a positive
solution (for other values of $b$ and $a$, a more detailed investigation
of the equation $f(t)=0$ may still yield solutions $t>0$, of course).
In both cases, the Ricci-flatness now follows from Proposition
$\ref{Prop: Sasaki-Ricci-Flat}$.

\qed

\begin{remark}
Notice that $t=1$ never gives a solution of
$f_1(t)=0$ that is, if $(M,g,\varphi,\xi,\eta)$ is 
Einstein then it is 
never $\nabla$-Einstein with respect to its characteristic connection,
in accordance with Proposition \ref{Prop:relation}  and the ensuing
discussion of normal forms in dimensions $5$ and $7$.
\end{remark}

\begin{remark}
 Observe that Theorem \ref{Theo: Einstein-Sasaki} leads to many 
examples of homogeneous $\nabla$-Ricci flat manifolds which are 
not flat, as opposed to the Riemannian case (for the standard case 
refer, for example, to \cite{Bess87i}). 
\end{remark}

We do not know of a similar general result for constructing explicitely
$\Ric^\nabla$-flat manifolds in \emph{even} dimensions.
Nevertheless, such manifolds exist, see Example \ref{CYT}.

%----------------------------------------------------------------------------
\subsubsection{$G_2$ manifolds with torsion}\label{section:G2}
%------------------------------------------------------------------------------
%
We now consider the class of 7-dimensional Riemannian manifolds
equipped with a $G_2$ structure. A $G_2$ structure can be seen as a
triple $(M,g,\omega)$ consisting of a 7-dimensional dimensional
manifold, a Riemannian metric, and a 3-form of general type at any
point. A $G_2$ $T$-manifold  --- $G_2$ manifold with (skew) torsion ---
is a manifold equipped with a $G_2$ structure such that there exists
a one-form $\theta$ such that $d\ast \omega = \theta \wedge
\ast\omega$; equivalently, these are the manifolds of
Fernandez-Gray type $\mathcal{X}_1\oplus\mathcal{X}_3\oplus \mathcal{X}_4$,
see \cite{FernandezGray82}. It admits  a unique connection with
totally skew symmetric torsion which preserves both the metric $g$
and the $3$-form $\omega$, called again the characteristic connection
(of the $G_2$ structure), see \cite{FriedIvan02}.

\begin{example}
%--------------
A $G_2$ manifold $(M,g,\omega)$ is called \emph{nearly parallel $G_2$}
if $d\omega = \lambda *\omega$, for some $0\neq \lambda\in\R$.
They coincide with the $G_2$ manifolds of type $\mathcal{X}_1$.
It is a well known fact that nearly parallel $G_2$ 
manifolds are Einstein\footnote{ 
This can be found in \cite{FKMS97}; 
the result is also implicitly contained in \cite[Thm. 13, p.120]{BFGK91}, 
since one checks that the assumptions of the theorem are  exactly describing
the existence of a nearly parallel $G_2$ structure. The result then
follows from the fact that every spin manifold with a Killing spinor
is Einstein, \cite{Frie80}. In that time, it was just not yet
fashionable to call $G_2$ structures by this name.}.  
Furthermore, a nearly parallel $G_2$ manifold is also
Einstein with parallel skew torsion,
$\Ric^\nabla = \frac{\lambda^2}{3} g$.
This can be easily deduced from the formulas in 
\cite[p.~318]{FriedIvan02} 
or from the spinorial argument in \cite{Agri06}, but 
it is also an immediate consequence of Proposition \ref{Prop:relation}.
\end{example}

\begin{example}
%--------------
There are $7$-dimensional cocalibrated $G_2$ manifolds $(M^7,g,\omega^3)$ 
with characteristic torsion $T$ such that
$$
\nabla T  =  0  , \quad d T =  0, \quad \delta T  =  0 , \quad 
\Ric^{\nabla}  =  0 , \quad \mathfrak{hol}(\nabla) \subset 
\mathfrak{u}(2)  \subset  \g_2  . 
$$
These $G_2$ manifolds have been described in \cite[Thm 5.2]{Frie07} (the
degenerate case where $ 2a +  c = 0$). $M^7$ is the
product $X^4 \times S^3$, where $X^4$ is any Ricci-flat K\"ahler manifold and
$S^3$ the round sphere. There is an underlying $S^1$-fibration, but
it does not induce a contact structure (the $1$-form $\eta\simeq e_7$ describing
the fibre satisfies $d\eta\wedge d\eta=0$, which cannot be); hence,
the example is not covered by the Tanno deformations described in 
Theorem \ref{Theo: Einstein-Sasaki}.
\end{example}

 $3$-Sasakian manifolds $M$ are Riemannian Einstein
spaces of positive scalar curvature carrying three compatible orthogonal 
Sasakian structures $(\xi_{\alpha} , 
\eta_{\alpha}, \varphi_{\alpha}), \alpha=1,2,3$ 
\cite{Ishihara,BoyerGalickiBook}.
The $7$-dimensional case is of particular interest because of its
relation to spin geometry and $G_2$ structures. The simply connected
homogeneous $3$-Sasakian $7$-manifolds were classified up to isometry in
\cite{Boyer&G&M94} and turn out to be
$\mathrm{Sp}(2)/\mathrm{Sp}(1)\cong S^7$ and  $\mathrm{SU}(3)/S^1$; 
by \cite{FriedrichKath90}, these are precisely the
compact simply connected 
$7$-dimensional spin manifolds with regular $3$-Sasaki structure.
Many non homogeneous examples are known \cite{BoyerGalickiBook}.

\begin{theorem}\label{thm:3-Sasaki-Einstein}
%--------------------------------------------
Every  $7$-dimensional $3$-Sasakian manifold carries three different
connections that turn it into an Einstein manifold with parallel skew
torsion; furthermore, it admits a deformation of the metric that
carries again an Einstein structure with parallel skew
torsion.

\end{theorem}

\proof 
%--------------
It is well-known that every $7$-dimensional $3$-Sasakian manifold admits
three linearly independent Killing spinors $\psi_1,\psi_2,\psi_3$ 
\cite{FriedrichKath90};
hence, each of these Killing spinors defines a nearly parallel $G_2$ structure
$\omega_i$ with a characteristic connection $\nabla^i$ which is Einstein
with parallel skew torsion by the previous example. As described in
\cite[Thm.~6.2]{AgriFrie10b}, these three connections are truly different,
and have the same constant $\lambda=-4$, hence they satisfy
$\Ric^\nabla = \frac{16}{3}g $.

Furthermore, any $7$-dimensional $3$-Sasakian metric can be deformed into
a $G_2$-structure in the following way (see \cite{FKMS97,Frie07,AgriFrie10b}
for more details). 
The vertical subbundle $\mathrm{T}^v \subset TM$ is spanned by 
$\xi_1, \xi_2, \xi_3$, its orthogonal complement is the horizontal 
subbundle $\mathrm{T}^h$. Fix a positive parameter $s > 0$ and
consider a new Riemannian metric $g^s$ defined by
$$
g^s(X  ,  Y) \ := \ g(X  ,  Y) \quad \mbox{if} \quad X,Y \in 
\mathrm{T}^h  , \quad 
g^s(X , Y) \ := \ s^2 \cdot g(X , Y) \quad \mbox{if} \quad X,Y \in 
\mathrm{T}^v .
$$
We rescale the $3$-forms 
$$F_1 =\eta_1 \wedge \eta_2 \wedge \eta_3  \quad \mbox{and} \quad
F_2 = \frac{1}{2}\big( \eta_1 \wedge d \eta_1 + \eta_2 \wedge 
d \eta_2  + \eta_3 \wedge d \eta_3 \big)  + 3  \eta_1 \wedge \eta_2
\wedge \eta_3 
$$
to obtain the new forms
$$
F_1^s \ := \ s^3 F_1  , \quad F_2^s \ := \ s F_2 , \quad \omega^s \ := \
F_1^s  +  F_2^s .
$$
One shows that
$(M^7, g^s, \omega^s)$ is a Riemannian $7$-manifold equipped with a
cocalibrated $G_2$-structure $\omega^s$, hence it admits a characteristic 
connection  $\nabla$ with skew torsion
$$
H_s \ = \ \left[ \frac{2}{s} - 10 s\right] (s\eta_1) \wedge (s \eta_2)
\wedge (s \eta_3) + 2 s \omega^s .
$$
The Ricci tensor of the characteristic connection of $\nabla$ (written as
an endomorphism) is given by the formula \cite{Frie07} 
$$
\Ric^\nabla \ = \ 12 \, (1 \, - \, s^2 ) \, 
\mathrm{Id}_{\mathrm{T}^h}  \, \oplus \, 
  16 \, (1 \, - \, 2 \, s^2) \, \mathrm{Id}_{\mathrm{T}^v}   .
$$
If $ s = 1$ (the $3$-Sasakian case), then $\mathrm{Ric}^{\nabla}$ vanishes
on the subbundle $\mathrm{T}^h$. For $s = 1/\sqrt{5}$, the Ricci tensor is
proportional to the metric, $\mathrm{Ric}^{\nabla}= (48/5) \,
\mathrm{Id}_{TM^7}$, so we obtain an Einstein structure with skew torsion
as claimed. As already observed in \cite[Thm.~5.4]{FKMS97}, it is nearly 
parallel $G_2$,
so it has parallel torsion. Finally, it was shown 
in \cite[Cor.~7.1]{AgriFrie10b}
that its corresponding Killing spinor is $\psi_0$, the \emph{canonical spinor}
introduced in the same paper. Since it is related to the
Killing spinors of the underlying $3$-Sasakian structure by
$\psi_i=\xi_i\cdot \psi_0$ \cite[Thm.~6.1]{AgriFrie10b}, 
we see that we obtained yet another different nearly parallel $G_2$
structure on the $3$-Sasaki manifold we started with.

\qed

%In the next paragraph, we will explicitely describe these Einstein structures
%with parallel skew torsion on the Aloff-Wallach manifold $N(1,1)$.

%----------------------------------------------------------------------------
\subsubsection{The Aloff-Wallach manifold 
$N(1,1)$}\label{AW}
%------------------------------------------------------------------------------
%

This is a computer-aided systematic search for metric connections with
skew torsion on the Aloff-Wallach manifold $\SU(3)/S^1$. It was our goal
to test how rare or common $\nabla$-Einstein structures are, apart from
the ones that we can predict theoretically. The main result is that, indeed,
additional $\nabla$-Einstein structures exist.

We use the computations available in \cite[p.109 ff]{BFGK91} and 
\cite[p.733 ff]{AgriFrie04}, which we then shall not reproduce here.
The manifold $N(1,1)$ is the homogeneous space
$\SU(3)/S^1$ where we are taking the embedding $S^1\rightarrow \SU(3)$ 
given by $\mathrm{e}^{i\theta} \mapsto 
\mathrm{diag}(\mathrm{e}^{i\theta},\mathrm{e}^{i\theta},\mathrm{e}^{-2i\theta})$.
The Lie algebra $\mathfrak{su}(3)$ splits into $\mathfrak{su}(3) 
= \mathfrak{m}+\R$, where $\R$ is the Lie algebra of $S^1$ given by 
the considered embedding. The space $\mathfrak{m}$ splits into 
$\mathfrak{m}_0\oplus \mathfrak{m}_1\oplus\mathfrak{m}_2\oplus
\mathfrak{m}_3$, where all $\mathfrak{m}_i$ are pairwise orthogonal 
with respect to the (negative of the) Killing form 
$B(X,Y):=-\mathrm{Re}(\mathrm{tr}\, XY)/2$ of $\mathfrak{su}(3)$.
The subspace
$\m_0$ is spanned by the matrix $L:=\mathrm{diag(3i, -3i,0)}$. 
Let $E_{ij}\, (i<j)$ 
be the matrix with $1$ at the place $(i,j)$ and zero 
elsewhere, and define $A_{ij}=E_{ij}-E_{ji},\,\tilde{A}_{ij}=i(E_{ij}+E_{ji})$.
Then  $\m_1:=\mathrm{Span}\{A_{12}, \tilde{A}_{12}\}$,  
$\m_2:=\mathrm{Span}\{A_{13}, \tilde{A}_{13}\}$ and 
$\m_3:=\mathrm{Span}\{A_{23}, \tilde{A}_{23}\}$. 
We consider 
the two-parameter family of metrics defined by the formula
$$
g_{s,y}:= \dfrac{1}{s^2}B|_{\mathfrak{m}_0}+B|_{\mathfrak{m}_1}
+\dfrac{1}{y}B|_{\mathfrak{m}_2}+ \dfrac{1}{y}B|_{\mathfrak{m}_2}.
$$  
This  is a subfamily
of the family considered in \cite[p.109 ff]{BFGK91}.
The isotropy representation $\Ad(\theta)$ leaves the vectors in 
$\mathfrak{m}_0$ and $\mathfrak{m}_1$ invariant,
and acts as a rotation by
$3\theta$ in the $\mathfrak{m}_2$-plane and in the $\mathfrak{m}_3$-plane. 
We use the standard realization of the
$8$-dimensional $\Spin(7)$-representation $\Delta_7$ as given in 
\cite[p.97]{BFGK91}, and denote by
$\psi_i,i=1,\ldots 8$ its basis. 
One then checks that $\psi_3,\psi_4,\psi_5$ and $\psi_6$
are fixed under the lift $\tilde{\Ad}(\theta)$ of the isotropy representation 
to $\Spin(7)$. Thus, they define constant sections in the spinor bundle
$\Sigma N(1,1)=\SU(3)\x_{\tilde{\Ad}}\Delta_7$. The metric defined by
$s=1,y=2$ is exactly the $3$-Sasakian metric on $N(1,1)$ (see the comments
in the previous section); it has three
Killing spinors with Killing number $1/2$ ($\psi_3,\psi_4,\psi_6$
in our notation). The metric defined by $s=1,y=2/5$ is the 
Einstein metric with  Killing spinor $\psi_5$ and with Killing number $-3/10$,
the well-known nearly parallel $G_2$ structure on $N(1,1)$
(see \cite[Thm 12, p.116]{BFGK91}).
It coincides with the nearly parallel $G_2$ structure constructed by
rescaling the underlying  $3$-Sasakian structure as described
in Theorem \ref{thm:3-Sasaki-Einstein}. 

In dimension $7$, any connection $\nabla$ with skew 
torsion $H$ admitting a parallel spinor field defines a $G_2$ structure 
of Fernandez-Gray type $\mathcal{X}_1\oplus\mathcal{X}_3\oplus \mathcal{X}_4$
on this manifold, and vice versa. This construction principle was used
in \cite{AgriFrie04} to define $G_2$ structures on $N(1,1)$ via their
parallel spinors. In particular, torsion forms $T_i$ 
depending on the metric parameters $s,y$ were given in 
Propositions 8.1--8.4 that admit $\psi_i$
as parallel spinors, $i=3,\ldots,6$.
Let $\nabla^i$ denote the connection with skew torsion $T_i$.
We can summarize our results as follows:

\begin{theorem}
%-----------
On the Aloff-Wallach manifold $N(1,1)$ with the family of metrics
$g_{s,y}$, the  connection $\nabla^i$, $i=3,\ldots,6$, defines
a $G_2$ structure that is  Einstein with skew torsion
precisely for the values stated in the following table:
$$
\begingroup
\renewcommand*{\arraystretch}{1.5}
\begin{array}{|l|l|l|l|l|}
%-----------------------
\hline
\text{Connection} & \text{Metric }g_{s,y} & s^\nabla/7 & \text{Comment} 
& \text{Riemannian Einstein ?}\\ \hline\hline
\nabla^3, \ \nabla^4 & s=1, y=2 & & \text{3-Sasakian} & \text{yes}\\ \hline
\nabla^5 & s=1, y=\frac{2}{5} & & \text{nearly parallel } G_2 & \text{yes}\\ \hline
\nabla^5 & s=1,y=\frac{30}{13} & \frac{7245}{1352} \cong 5.3587 & \text{new} &  \text{no}\\ \hline
\nabla^6 & s=1, y=2 & & \text{3-Sasakian} & \text{yes}\\ \hline
\nabla^6 & s\cong 0.97833, & \cong 1.676989544 & \text{new} &  \text{no}\\ 
& y\cong 0.34935 &&&\\ \hline
\end{array}
\endgroup
$$
\end{theorem}

Hence, we were able to construct two Einstein structures with skew
torsion on $N(1,1)$ that go beyond the metrics predicted for theoretic
reasons. It will be an interesting topic for further research to investigate
their detailed geometrical properties.

\begin{remark}
%--------------
In  \cite{AgriFrie04}, one can also find a construction of $3$-forms
that make various linear combinations of the spinors $\psi_3,\ldots,\psi_6$
parallel. The algebraic systems of equations for the $\nabla$-Einstein 
condition are difficult to control, but tests with different parameters 
lead us to conjecture that for arbitrary $a,b$ $(ab \neq 0)$ and
$a,b,c$ $(abc \neq 0)$ we obtain no new solutions. 

In a same vein, one can ask for $\nabla$-Einstein structures
for other embeddings of $S^1$ into $\SU(3)$, i.\,e.~the general
Aloff-Wallach manifold  $N(k,l)$, with $\ k,l\in\mathbb{N}$ coprime. 
For these, already the
Riemannian Einstein metrics cannot be described explicitly,
but only as solutions of a complicated system of equations
depending on $k,l$, and the metric parameters $s,y$.
Some tests for different value of $k$ and $l$ showed that
further $\nabla$-Einstein structures do exist, but it seems hopeless
to discuss the resulting system of equations in a reasonably
general way for arbitrary $k$ and $l$.
\end{remark}

%----------------------------------------------------------------------------
\appendix\section{Appendix: Equivalent formulations of parallel torsion}
%----------------------------------------------------------------------------
Consider the one parameter family of connections with skew torsion
given by
$$\nabla^s_X Y = \nabla^g_X Y + 2s H(X,Y), \quad s\in \R.$$
We will investigate the relation between parallel torsion and 
the vanishing of the first Bianchi identity.
For the connection $\nabla^s$ the following curvature identity holds
$$\begin{array}{lcl} R^s(X,Y,Z,W) & = & R^g(X,Y,Z,W) + 4s^2 g(H(X,Y),H(Z,W))+ 4s^2 \sigma_H (X,Y,Z,W)\medskip\\
& &+ 2s \nabla^s_X H (Y,Z,W) - 2s \nabla^s_Y H (X,Z,W)\end{array}$$
and using the fact that 
$$\nabla^s_X H(Y,Z,W) = \nabla^g_X H (Y,Z,W) - 2s\sigma_H(X,Y,Z,W)$$
it is easy to check that the first Bianchi identiy reads as
$$\stackrel{XYZ}{\sigma}R^s(X,Y,Z,W) = -8 s^2 \sigma_H(X,Y,Z,W) + 4 s (dH(X,Y,Z,W) + \nabla^g_W H(X,Y,Z)).$$

\begin{proposition}\label{Prop: Bianchi}
For any real parameter $s_0\neq 0$, the following conditions are equivalent:
 \begin{enumerate}
  \item $\nabla^{s_0}H =0$ ;
  \item  $\stackrel{XYZ}{\sigma}R^{3 s_0}(X,Y,Z,W)=0$;
  \item $dH=8s_0\sigma_H$.
  \end{enumerate}
\end{proposition}
\proof 
%-------------
Suppose $(1)$ holds, i.\,e.~there exists $s_0\neq 0$ such that $\nabla^{s_0}H = 0$. 
Then $\nabla^g H =  2 s_0 \sigma_H$, which implies that $dH = 8 s_0
\sigma_H$. The Bianchi identity then reduces to
$$
\stackrel{XYZ}{\sigma}R^{s}(X,Y,Z,W) = -8 s (s-3s_0)\sigma_H(X,Y,Z,W).
$$ 
Hence, we have $(2)$, the vanishing of the first Bianchi identity for 
$s = 3 s_0$. 

Consider now condition $(2)$.
The vanishing of the first Bianchi identity for some $s\neq 0$ implies that 
$$
\nabla^g_W H(X,Y,Z) = 2 s \sigma_H (X,Y,Z,W) - dH(X,Y,Z,W)
$$
and from this expression we see that $\nabla^g H$ is a totally
antisymmetric tensor, and is therefore equal to $\frac{1}{4}dH$.
Then 
$$
dH(X,Y,Z,W) + \frac{1}{4} dH(W,X,Y,Z) = 2s\sigma_H(X,Y,Z,W)
$$ and simplifying we conclude that
$dH = \frac{8}{3} s\, \sigma_H,$, that is, $(3)$ holds.
But then the first equation can be rewritten
$\nabla^g H = \frac{2}{3} s \sigma_H$, which is  equivalent to 
$\nabla^{s/3} H =0$, i.\,e.~we proved $(1)$.

It remains to deduce $(1)$ from $(3)$. For this, consider the general
identity (\cite{FriedIvan02}, \cite[Cor. A.1.]{Agri06})
$$
dH(X,Y,Z,W) + \nabla^s_W H(X,Y,Z) - 8s\sigma_H(X,Y,Z,W)
\ =\ \stackrel{XYZ}{\sigma} [\nabla^s_X H(Y,Z,W)] ,
$$
which holds under the assumption that $s\neq 0$.
Both sides are tensorial quantities, but the left hand side is symmetric
in $X,Y,$ and $Z$, while the right hand side is antisymmetric in 
$X,Y,$ and $Z$, hence they
can only be equal if they vanish, 
i.\,e.~$\stackrel{XYZ}{\sigma} [\nabla^s_X H(Y,Z,W)]=0$ and
\begin{equation}\label{nice}
dH(X,Y,Z,W) + \nabla^s_W H(X,Y,Z) - 8s\sigma_H(X,Y,Z,W)
\ =\ 0. 
\end{equation}
Assuming $(3)$, this identity
is reduced to  $\nabla^s_W H(X,Y,Z)=0$, which is exactly
condition $(1)$.
 
\qed

Proposition \ref{Prop: Bianchi} is clearly wrong for the Levi-Civita
connection ($s_0=0$), showing the non-triviality of the result.
In dimension $4$, the $4$-form $\sigma_H$ vanishes for purely
algebraic reasons; the theorem stays correct and says basically
that $dH=0$ is equivalent to $\nabla H=0$. Observe that in the case of
non-parallel torsion, identity \ref{nice} still holds ($s\neq 0$) and is, in our
opinion, quite remarquable: it generalizes in a straight forward way
conditions (1) and (3), and it proves that any two of the three quantities
$\nabla H, dH$, and $\sigma_H$ determines the third.
%

%\bibliographystyle{amsalpha}
%\bibliography{bibliography}

\providecommand{\bysame}{\leavevmode\hbox to3em{\hrulefill}\thinspace}
\providecommand{\MR}{\relax\ifhmode\unskip\space\fi MR }
\providecommand{\MRhref}[2]{%
  \href{http://www.ams.org/mathscinet-getitem?mr=#1}{#2}
}
\providecommand{\href}[2]{#2}

\newpage

\begin{center}
\textsc{\textbf{\Large Erratum}} 

\vspace*{0.5 cm}

{\large Ilka Agricola \& Ana Cristina Ferreira} 

\vspace*{0.5 cm}

{\large June 29, 2022.}

\vspace*{1cm}

{\bf Abstract:} We correct here a wrong argument which appeared in example 2.14.

\vspace*{1cm}

\end{center}

\section*{3-forms of Einstein type}
Let $(M,g)$ be a Riemannian manifold and $H\in\Lambda^3(M)$.
The  (unique) metric connection with skew torsion $H$ is
given by$$
\nabla_X Y = \nabla^g_X Y +\frac{1}{2}H(X,Y,-).
$$
We denote quantities refering to the Levi-Civita connection with
an upper index $g$, while quantities associated with the new connection
will have an upper index $\nabla$.  

In subsection \ref{subsec: notations}, we introduced the tensor
\begin{equation*}\label{diff-tensor-T}
S(X,Y) :=  \sum_{i=1}^n g(H(e_i,X),H(e_i,Y)) 
= \sum_{i,j=1}^n H(e_i,X,e_j)H(e_i,Y,e_j)
\end{equation*}
where $(e_1, \cdots,e_n)$ is an orthonormal frame of $TM$. This tensor measures the (symmetric part of the) difference between
the Riemannian and the $\nabla$-curvature.  Indeed, we have that 

$$\mathrm{Ric}^\nabla(X,Y)  =  \Ric^g(X,Y)-\frac{1}{4} 
 S(X,Y) - \frac{1}{2} \delta H (X,Y).$$

\begin{definition}
%------------------
On a Riemannian manifold $(M,g)$, a $3$-form $H$ will be called `of 
Einstein type' if the difference tensor $S(X,Y) := \sum_i g(H(e_i,X),
H(e_i,Y))$ is proportional to the metric $g$.
\end{definition}

We proved the following combinatorial criterion.

\begin{proposition}\label{Prop:relation}
%----------------------------------------
Let $(M,g)$ be a Riemannian manifold, $H$ a $3$-form 
written in a local orthonormal frame $e^1,\dots,e^n$ of $T^*M$,
$
H = \sum_{ijk} H_{ijk}\, e^i \wedge e^j \wedge e^k.
$ 
Then  $H$ is of Einstein type if it satisfies the following conditions:
\begin{enumerate}
\item no term of the form $H_{ija} e^i\wedge e^j \wedge e^a + H_{ijb} e^i\wedge e^j \wedge
e^b$ with $a\neq b$ occurs;
\item if $i$ and $j$ are two indices in $\{1,\dots, n\}$ then the
number of occurrences of $i$ and $j$ in $H$ coincides;
\item if $\{i,j,k\}$ and $\{a,b,c\}$ are two sets of indices then $H_{ijk}^2 =
H_{abc}^2$.
\end{enumerate}
\end{proposition}

\begin{remark} 
%--------------
Proposition $\ref{Prop:relation}$ yields an easy procedure for
producing further examples of $\nabla$-Einstein metrics with non
constant scalar curvature 
for all manifolds that are parallelizable and carry an Einstein metric 
(for example $S^7$ or compact semi-simple Lie groups).
\end{remark}

\begin{figure}[t]
\begin{center}
\psfrag{I}{$e_{123}$}
\psfrag{II}{$e_3(e_{15}+e_{24})$}
\psfrag{III}{$e_{126}+e_{135}+e_{234}$}
\psfrag{IV}{\!\!\!${}^{*}\ e_{123}+e_{456}$}
\psfrag{V}{\!\!\!\!\!${}^{*}\ e_{135}+e_{146}+e_{256}+e_{234}$}
\psfrag{6}{VI. $e_1(e_{23}+e_{45})+e_{267}$}
\psfrag{7}{VII. $e_{123}+e_{456}+e_7(e_2+e_5)(e_3+e_6)$}
\psfrag{8}{VIII. $e_1(e_{23}+e_{45}+e_{67})$}
\psfrag{9}{IX.$^{**}$ VIII + $e_{246}$,}
\psfrag{10}{X. VII + $e_{147}$}
\psfrag{11}{XI. VIII + $e_2(e_{46}- e_{57})$}
\psfrag{12}{XII. $e_1(e_{45}+e_{67})+e_2(e_{46}-e_{57})+ e_3(e_{47}+e_{56})$,}
\psfrag{13}{XIII.$^{**}$ XII $- e_{123}$}
\includegraphics[width=8cm]{schiouten-ergaenzt.eps}
\end{center}
\caption{Normal forms of $3$-forms up to dimension $7$ (see Example 
\ref{exa:normal-forms}); the picture is
taken from the original article of J.~Schouten from 1931 \cite{Schouten31}
while for $n=7$ we followed \cite{Westwick81}. Any three numbers $a,b,c$ linked
by lines represent a summand $e_a\wedge e_b\wedge e_c$.}
\label{normal-forms}
\end{figure}

\begin{example}\label{exa:normal-forms}
%--------------------------------------
In 1931, J.~Schouten described the normal forms of 
$3$-forms up to dimension $7$ \cite{Schouten31}, i.\,e.~representatives of the 
$\mathrm{GL}(n,\R)$-orbits for $n\leq 7$, see Figure \ref{normal-forms} 
(see \cite{Westwick81} 
for a modern account of the real classification). 

We argued \emph{wrongly} that one checks
by a direct computation that only the following $3$-forms 
are of Einstein type in the given dimensions:
\begin{enumerate}
\item Type I. in dimension $3$,
\item Types IV. and V. in dimension $6$,
\item Type XIII. in dimension $7$.
\end{enumerate}
and that, in particular, Riemannian Einstein manifolds $(M,g)$ will never be
$\nabla$-Einstein in dimensions $4$ and $5$. 
\end{example}

Clearly,  such a claim cannot be made since, for each dimension, $\mathrm{O}(n)$ normal forms should  have been considered instead of $\mathrm{GL}(n,\mathbb{R})$ ones.  There is a priori no guarantee that the $\mathrm{GL}(n,\mathrm{R})$ are written with respect to an orthornormal frame or that other normal forms do not arise.

However, by writting the three-form $H$ in full generality as $H= \sum_{i,j,k =1} H_{ijk} e^i\wedge e^j \wedge e^k$ and using the combinatorial criterion, it is possible to prove that the claimed results hold in dimensions 4, 5 and 6. We did not check in dimension 7 since the space of 3-forms has dimension 35 and becomes too large.

We finish this note by remarking that this incorrect example does not impact the rest of the manuscript.

\bigskip

\medskip\noindent {\bf Acknowledgements.} 
Both authors thank Jorge Lauret (C\'ordoba) for his interest in our work and for pointing out the mistake in question.

\end{document}